\documentclass[12pt]{article}
\usepackage[final]{epsfig}
\usepackage{graphics}
\usepackage{amsmath}
\usepackage{amsfonts}
\usepackage{latexsym}
\usepackage{amssymb}
\usepackage{graphicx}
\usepackage{url}

\newtheorem{lemma}{Lemma}[section]
\newtheorem{proposition}[lemma]{Proposition}
\newtheorem{remark}[lemma]{Remark}

\newtheorem{theorem}{Theorem}

\newtheorem{corollary}[lemma]{Corollary}

\newcommand{\g}{{\gamma}}
\newcommand{\G}{{\Gamma}}
\newcommand{\proofend}{$\Box$\bigskip} 

\newcommand{\C}{{\mathbb C}}

\newcommand{\R}{{\mathbb R}}

\newcommand{\RP}{{\mathbb {RP}}}  

\def\proof{\paragraph{Proof.}}

\begin{document}

\title{On the discrete bicycle transformation}

\author{S. Tabachnikov\footnote{
Department of Mathematics,
Pennsylvania State University, 
University Park, PA 16802, USA, 
tabachni@math.psu.edu} \and E. Tsukerman\footnote{Stanford University, emantsuk@stanford.edu}
}

\date{}
\maketitle

\section{Introduction} \label{intro}

The motivation for this paper comes from the study of a simple model of bicycle motion. The bicycle is modeled 
as an oriented segment in the plane of fixed length $\ell$, the wheelbase of the bicycle.  The motion is constrained so that the segment is always tangent to the path of the rear wheel; this non-holonomic constraint is due to the fact that the rear wheel is fixed on the frame, whereas the front wheel can steer. See \cite{Fi,FLT,LT,Ta1} and the references therein.

\begin{figure}[hbtp]
\centering
\includegraphics[width=2.4in]{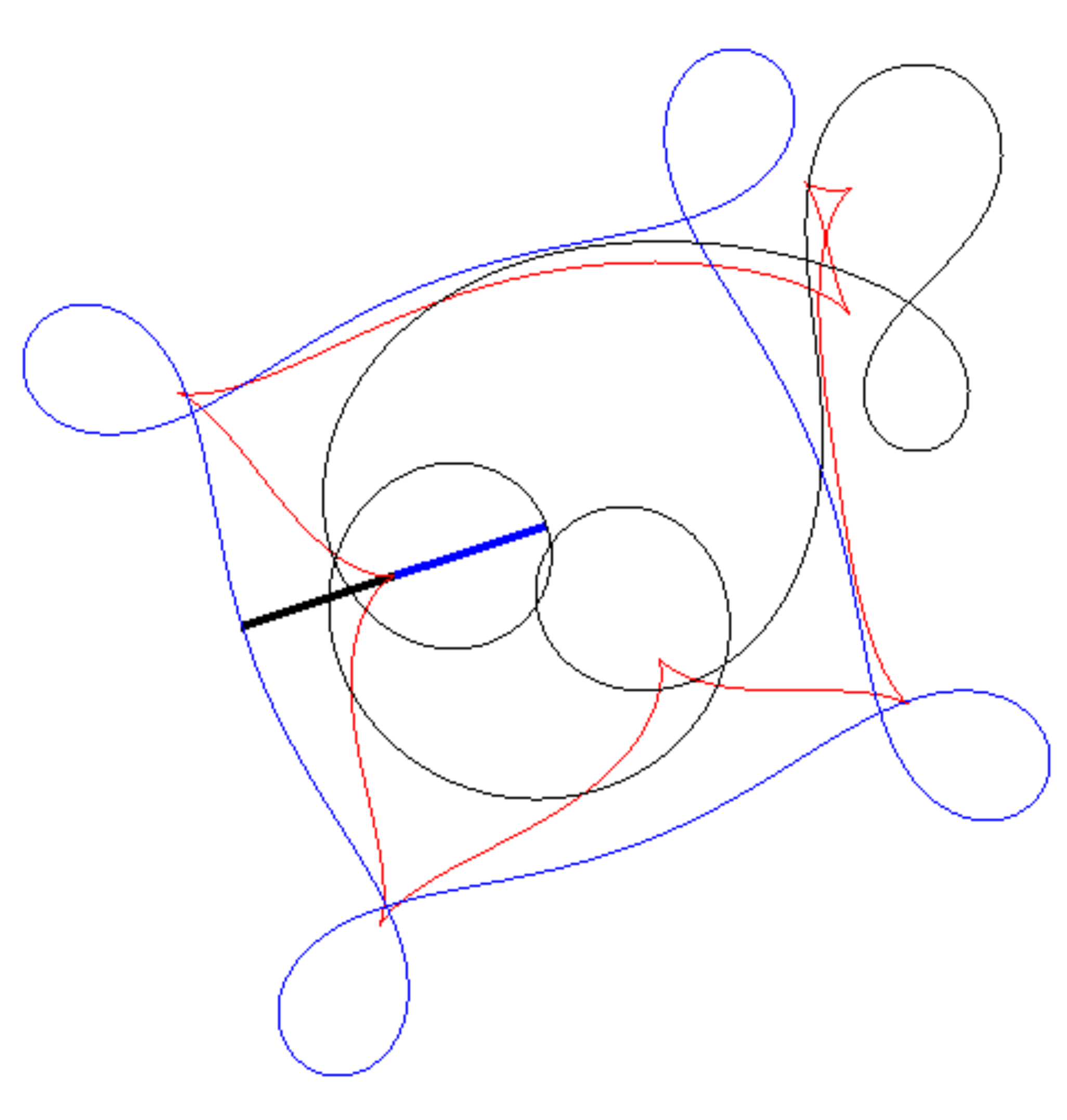}
\caption{Bicycle correspondence. The cusped curve is the rear track, the two smooth curves are front tracks in the bicycle correspondence (figure courtesy of R. Perline).}
\label{Darboux}
\end{figure}

If the rear wheel path $\g$ is prescribed, and the direction of motion is chosen, the front wheel path $\G$ is  constructed by drawing the tangent segments of length $\ell$ to $\g$. Note that the rear track may have cusp: they occur when the steering angle equals $90^{\circ}$. Changing the direction of motion to the opposite yields another front track, say, $\G'$. We say that the curves $\G$ and $\G'$ are in the {\it bicycle correspondence}.\footnote{One can also call this Darboux or B\"acklund transformation, but we shall use the ``bicycle" terminology.} See Figure \ref{Darboux}.

If the front wheel path $\G$ is prescribed then the rear wheel follows a constant-distance pursuit curve, and its trajectory is uniquely determined, once the initial position of the bicycle is chosen. A monodromy map $M_{\G,\ell}$ arises that assigns to every initial position of the bicycle its terminal position. If $\G$ is a closed curve then $M_{\G,\ell}$ is a self-map of a circle of radius $\ell$, uniquely defined up to conjugation. The bicycle monodromy $M_{\G,\ell}$ is  a M\"obius transformation \cite{Fo,FLT,LT}.

All of the above can be extended to the motion of a segment in higher dimensional Euclidean spaces and even Riemannian manifolds (see \cite{HPZ} for elliptic and hyperbolic planes). In the forthcoming paper \cite{Ta2}, we shall discuss Liouville integrability of the 
bicycle transformation in dimensions 2 and 3.

\begin{figure}[hbtp]
\centering
\includegraphics[height=1.5in]{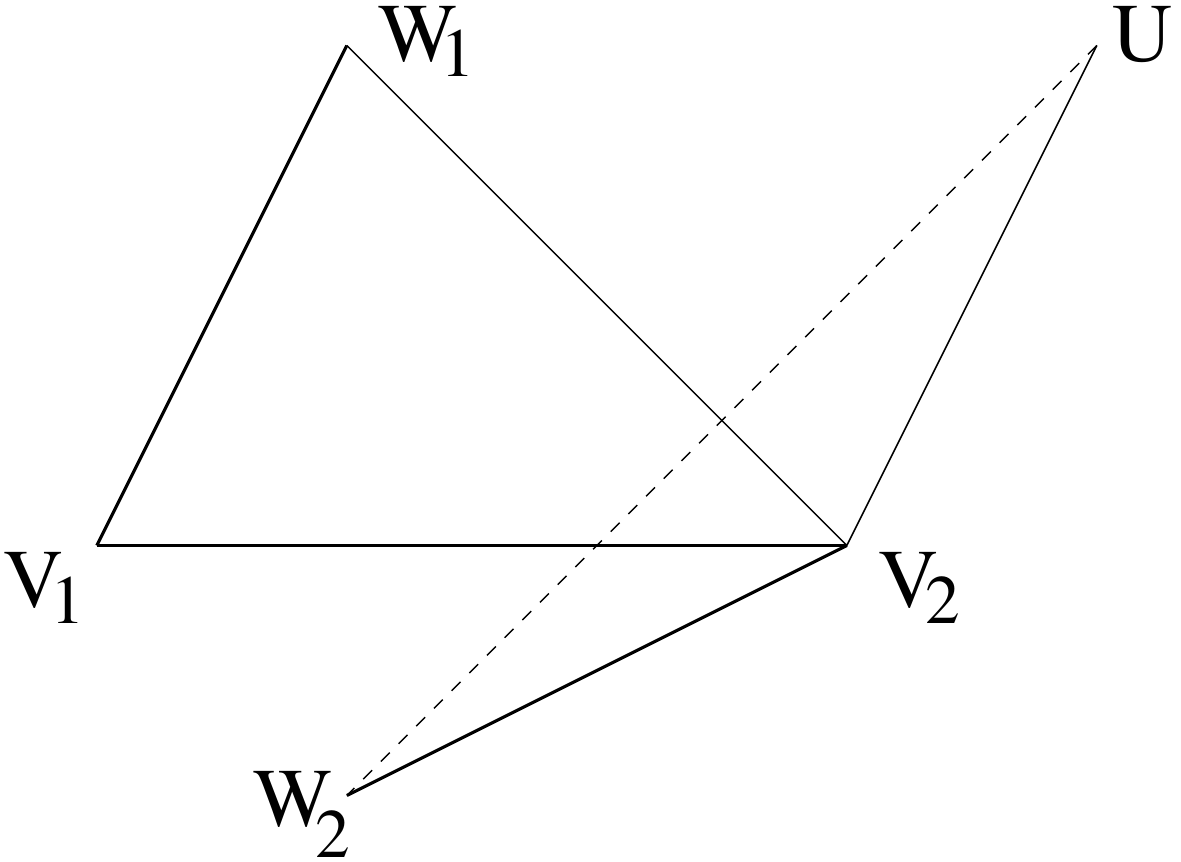}
\caption{Discrete bicycle correspondence}
\label{constr}
\end{figure}

In this paper, following \cite{Ho,PSW}, we study a discrete version of the bicycle correspondence. Let $V=(V_1, V_2, \dots)$ be a polygon in $\R^n$, and let $V_1 W_1$ be a seed segment of length $\ell$ (so now $\ell$ is twice the lenght of the bicycle frame). The next point $W_2$ is constructed in the plane spanned by $V_1, V_2, W_1$ as follows: one parallel translates point $W_1$ along the vector $V_1 V_2$ to point $U$, and then reflects point $U$ in the line $W_1 V_2$ to obtain a new point $W_2$. In other words, the plane quadrilateral $V_1 V_2 W_1 W_2$ is an isosceles trapezoid with $|V_1 V_2|=|W_1 W_2|$ and $|V_1 W_1|=|V_2 W_2|=\ell$, see Figure \ref{constr}. Once the point $W_2$ is constructed, one continues the process, shifting the index by one, etc.\footnote{The definition in \cite{Ho,PSW}, given in 3-dimensional case,  involves another, twist, parameter.} 

We call the above described correspondence between polygons $V$ and $W$ the {\it discrete bicycle correspondence} and denote it by ${\cal B}_{\ell} (V,W)$. In the continuous limit, the polygons $V$ and $W$ become  the front tire tracks $\G$ and $\G'$, and  the discrete bicycle correspondence becomes the above described bicycle correspondence between smooth curves. 

Our ultimate goal is to establish Liouville integrability of the discrete bicycle correspondence and to describe its dynamics in detail. In this paper, we make steps in this direction. Let us list basic properties of the discrete bicycle correspondence. 

Let $V$ be a closed $k$-gon in $\R^n$ (that is, $V_{i+k}=V_i$ for all $i$). The polygon $W$ is not necessarily closed, and the discrete bicycle monodromy $M_{V,\ell}$ arises, similarly to the continuous case. 

\begin{theorem} \label{monod}
The monodromy $M_{V,\ell}:S^{n-1} \to S^{n-1}$ is a M\"obius transformation of the sphere of radius $\ell$.
\end{theorem}

Thus, fixed points of the monodromy $M_{V,\ell}$ correspond to closed polygons $W$ in the  discrete bicycle correspondence with $V$. 

\begin{theorem} \label{conj}
Let $V$ and $W$ be closed polygons in $\R^n$ in the  discrete bicycle correspondence. Then, for every $\lambda$, the monodromies $M_{V,\lambda}$ and $M_{W,\lambda}$ are conjugated to each other.
\end{theorem}

Theorem \ref{conj} implies that the invariants of the conjugacy class of the monodromy, viewed as functions of the ``spectral parameter" $\lambda$, are integrals of the discrete bicycle correspondence. We refer to them as the monodromy integrals.

The next theorem states that the  discrete bicycle correspondences with different length parameters commute with each other (``Bianchi permutability"). Recall that we write ${\cal B}_{\ell} (V,W)$ to indicate that  polygons $V$ and $W$ are in the discrete bicycle correspondence with the length parameter $\ell$.

\begin{theorem} \label{comm}
Let $V, W, S$ be closed $k$-gons in $\R^n$ such that ${\cal B}_{\ell} (V,W)$ and ${\cal B}_{\lambda} (V,S)$ hold. Then there exists a closed polygon $T$ such that ${\cal B}_{\ell} (S,T)$ and ${\cal B}_{\lambda} (W,T)$ hold.
\end{theorem}

In the case of 3-dimensional space, Theorems \ref{monod}-\ref{comm} are not new:  in \cite{PSW}, they are proved using quaternions. We give different proofs in Section \ref{basic}. 

V. Adler \cite{Ad1,Ad2} studied complete integrability of a correspondence on the space of polygons in Euclidean space called the {\it recutting of polygons}. The recutting $R_i$ of polygon $V$ at $i$th vertex is the reflection of $V_i$ in the perpendicular bisector hyperplane of the segment $V_{i-1} V_{i+1}$ . Recuttings of $k$-gons form a group with generators $R_i,\ i=1,\dots, k$ and the relations
$$
R_i^2=1,\ R_i R_j =R_j R_i\ {\rm for}\ |i-j|\ge 2,\ {\rm and}\ R_i R_{i+1} R_i = R_{i+1} R_i R_{i+1},
$$
where the indices are understood cyclically.

The recutting is closely related to the discrete bicycle correspondence. In Section \ref{int}, we show that certain integrals of the recutting, discovered by Adler,  are integrals of the discrete bicycle correspondence. To do so, we construct a discrete analog of the rear track trajectory, a chain of mutually tangent spheres.

We also have the following result relating the discrete bicycle correspondence and the recutting.

\begin{theorem} \label{rec}
1) The monodromy is preserved by the recutting. In particular, the monodromy integrals are also integrals of the recutting.\\
2) The discrete bicycle correspondence commutes with the recutting.
\end{theorem}

\begin{figure}[hbtp]
\centering
\includegraphics[width=2in]{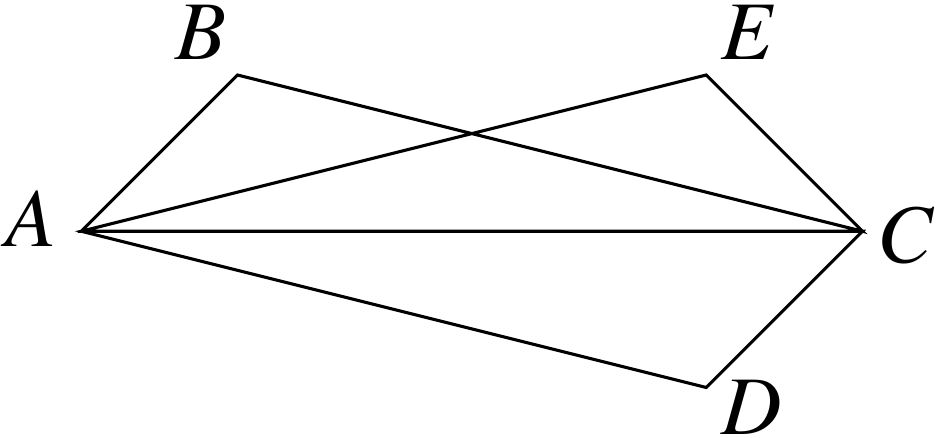}
\caption{The parallelogram $ABCD$ and the kite $AECD$ have the same monodromy}
\label{kite}
\end{figure}

To illustrate the first claim of Theorem \ref{rec}, a parallelogram and the corresponding kite have the same monodromy, see Figure \ref{kite}.

Theorems \ref{monod}-\ref{rec} are proved in Section \ref{basic}.

Consider the low-dimensional situation. If the dimension equals $2$ then the discrete bicycle monodromy belongs to $SL(2,\R)$. Then one has the trichotomy: $M_{V,\ell}$ may be elliptic, parabolic, and hyperbolic. In the last case, $M_{V,\ell}$ has two fixed points, and one can choose one (say, the attracting one) to construct a closed polygon $W$ in the discrete bicycle correspondence with $V$ (with length parameter $\ell$). According to Theorem \ref{conj}, $M_{W,\ell}$ is again hyperbolic, and one may iterate the construction by choosing the other fixed point of $M_{W,\ell}$ (otherwise, one gets back to $V$). Thus,  the discrete bicycle correspondence becomes a map on polygons, and we write ${\cal T}_{\ell} (V) = W$. 

In dimension three,  the discrete bicycle monodromy belongs to $SL(2,\C)$. If the monodromy is not the identity, it has two fixed points (perhaps, coinciding), and once again, one can consider the discrete bicycle correspondence as a mapping of the space of polygons in $\R^3$. 

In Sections \ref{plane} and \ref{quad}, we study the discrete bicycle transformation on plane polygons. We prove that the discrete bicycle transformation is defined on convex cyclic polygons only if the length parameter does not exceed the diameter of the circumcircle, and in this case, the transformation is a rotation about the circumcenter. We also compute the eigenvalues of the discrete bicycle monodromy and derive a criterion for the monodromy to be parabolic. 

In Section \ref{quad}, we  give a complete description of the dynamics of the discrete bicycle transformation on plane quadrilaterals. As an application, we classify the so-called bicycle $(4k,k)$-gons (see Section \ref{quad} for definition).

\section{Proofs of basic properties} \label{basic} 

\noindent{\bf Proof of Theorem \ref{monod}.}
Recall that the Mob\"ius group $O(n,1)$ consists of linear isometries of the pseudo-Euclidean space $\R^{n,1}$, and it acts projectively on $S^{n-1}$, the spherization of the null cone; it is also the group of isometries of $n$-dimensional hyperbolic space (in the hyperboloid model).

Let $M$ be the monodromy along segment $V_1 V_2$ in Figure \ref{constr}. We need to show that $M \in O(n,1)$.

Let $u,v$ and $x$ be the unit vectors along $V_1 W_1, V_2 W_2$ and $V_1 V_2$, respectively, and let $|V_1 V_2|=a$. The reflection of vector $u$ in vector $\xi$ is given by the formula
$$
v=\frac{2 u\cdot \xi}{|\xi|^2}\xi - u.
$$
Applying this to $\xi=ax-\ell u$, we obtain
\begin{equation} \label{refl}
v=\frac{u+\frac{2a^2 (x\cdot u)}{\ell^2-a^2}x-\frac{2a\ell}{\ell^2-a^2}x}{\frac{\ell^2+a^2}{\ell^2-a^2}-\frac{2a\ell (x\cdot u)}{\ell^2-a^2}}.
\end{equation}

On the other hand, a matrix from $O(n,1)$ has the form
$$
\left(\begin{array}{cc}
A&\xi\\
\eta^t&\lambda
\end{array}\right) 
$$
where $A$ is an $n\times n$ matrix, $\xi$ and $\eta$ are $n$-vectors, and the following relations hold:
$$
A^t A=E+\eta\otimes\eta^t,\ A^t (\xi)=\lambda \eta,\ \xi\cdot \xi=\lambda^2-1,
$$
where $E$ is the unit matrix, and $\eta\otimes\eta^t$ is the rank one matrix obtained by multiplying a column and a row vectors. The projective action of such a matrix is given by the formula:
\begin{equation} \label{proj}
u\mapsto \frac{A(u)+\xi}{\eta\cdot u +\lambda}.
\end{equation}

We observe that (\ref{refl}) has the form (\ref{proj}) with 
$$
A=E+\frac{2a^2}{\ell^2-a^2} x\otimes x,\ \xi=\eta=-\frac{2a\ell}{\ell^2-a^2} x,\ \lambda=\frac{\ell^2+a^2}{\ell^2-a^2},
$$
which completes the proof.
\proofend

 In dimension two, one identifies the unit circle  with the real projective line via stereographich projection from point $(-1,0)$. Then M\"obius transformations become fractional-linear. 
If $\alpha$ is the angular coordinate on $S^1$ then $x=\tan (\alpha/2)$ is the respective affine coordinate on $\RP^1$. In Figure \ref{constr}, assume that $V_1 V_2$ is horizontal, the direction of $V_1 W_1$ is $\alpha$ and that of $V_2 W_2$ is $\beta$. If $x=\tan (\alpha/2)$ and $y=\tan (\beta/2)$ then the monodromy is given by the formula
$$
y=\frac{\ell+a}{\ell-a}x,
$$
or
$$
M_{\ell}=
\left(\begin{array}{cc}
\ell+a&0\\
0&\ell-a
\end{array}\right). 
$$
In general, if the direction of $V_1 V_2$ is $\phi$ then
\begin{equation} \label{mono}
M_{\ell}=
\left(\begin{array}{cc}
\ell+a\cos \phi&-a \sin\phi\\
-a \sin\phi&\ell-a\cos\phi
\end{array}\right). 
\end{equation}

\begin{figure}[hbtp]
\centering
\includegraphics[height=.8in]{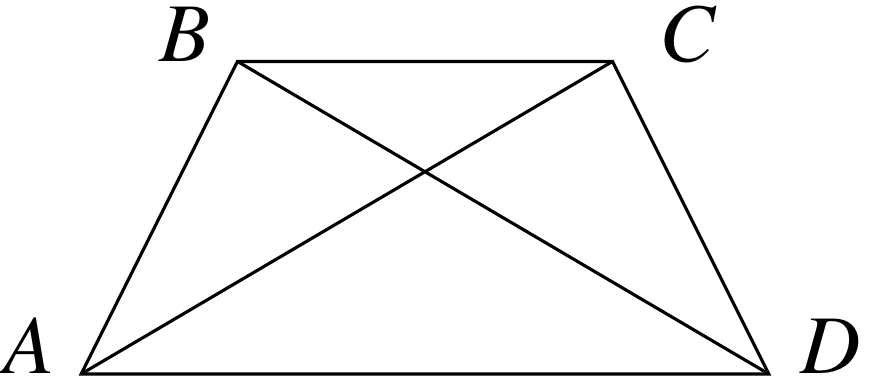}
\caption{A Darboux butterfly}
\label{trap}
\end{figure}

Now we prove a  property of isosceles trapezoids that is fundamental for what follows. Let $ABCD$ be a plane isosceles trapezoid, see Figure \ref{trap}. We call the closed quadrilateral $ABDC$, made of the lateral sides and diagonals of a trapezoid, a {\it Darboux butterfly}.

\begin{lemma}[Butterfly Lemma] \label{butterfly}
The monodromy (with any length parameter $\ell$) along a Darboux butterfly is the identity. Conversely, if the monodromy along a closed quadrilateral is the identity for some value of $\ell$ then the quadrilateral is a Darboux butterfly.
\end{lemma} 

\proof
The first statement of the lemma  is 3-dimensional: if $w$ is a test vector at vertex $A$ then the respective vectors at all other vertices (the ``transports" of $w$ along the quadrilateral) belong to the 3-dimensional space, spanned by the plane of the trapezoid and the vector $w$.

In fact, it suffices to consider the case when $w$ is in the plane of the trapezoid. Indeed, in dimension three, the monodromy is considered as an orientation preserving isometry of hyperbolic space acting on the sphere at infinity. If such an isometry has more than two fixed points then it is the identity.

In dimension two, we shall prove that the monodromy along the polygonal path $ABD$ equals the monodromy along the path $ACD$ if and only if $ABDC$ is a Darboux butterfly. Without loss of generality, assume that $AD$ is horizontal. Let $a,b,c,d$ be the length of the segments $AB, BD, AC, CD$, and let $\alpha,\beta,\gamma,\delta$ be the angles made with the positive horizontal axis. Let $|AD|=g$. 

The product of the  matrices from equation (\ref{mono}) is
$$
\left(\begin{array}{cc}
\ell-b\cos\beta & -b\sin\beta\\
-b\sin\beta & \ell+b\cos\beta
\end{array}\right)\left(\begin{array}{cc}
\ell-a\cos\alpha & -a\sin\alpha\\
-a\sin\alpha & \ell+a\cos\alpha
\end{array}\right),
$$
so we have the monodromy
\begin{equation} \label{matr}
M(a,b,\alpha,\beta)=\left(\begin{array}{cc}
\ell^{2}-lg+ab\cos(\alpha-\beta) & -ab\sin(\alpha-\beta)\\
ab\sin(\alpha-\beta) & \ell^{2}+lg+ab\cos(\alpha-\beta)
\end{array}\right).
\end{equation}
For equality to hold, we must have
$$
M(a,b,\alpha,\beta)=k(\ell)M(c,d,\gamma,\delta)
$$
for some constant $k(\ell)$ dependent only on $\ell$. Therefore 
$$
\frac{\ell^{2}-\ell g+ab\cos(\alpha-\beta)}{\ell^{2}-\ell g+cd\cos(\gamma-\delta)}=\frac{ab\sin(\alpha-\beta)}{cd\sin(\gamma-\delta)}=\frac{\ell^{2}+\ell g+ab\cos(\alpha-\beta)}{\ell^{2}+\ell g+cd\cos(\gamma-\delta)}.
$$
Set $X=\ell^{2}+ab\cos(\alpha-\beta)$ and $Y=\ell^{2}+cd\cos(\gamma-\delta)$.
Then 
$$
\frac{X-\ell g}{Y-\ell g}=\frac{X+\ell g}{Y+\ell g},
$$
hence $X=Y$ and
\begin{equation} \label{eqs}
ab\cos(\alpha-\beta)=cd\cos(\gamma-\delta), \ \ ab\sin(\alpha-\beta)=cd\sin(\gamma-\delta).
\end{equation}

The second equation (\ref{eqs}) implies that the signed
area of triangle $ABD$ is equal to that of triangle $ACD$,
so that the quadrilateral $ABDC$ has a total signed area of zero. 
It also follows that
$\tan(\alpha-\beta)=\tan(\gamma-\delta)$,
so that $\alpha-\beta=\gamma-\delta$ or $\alpha-\beta=\gamma-\delta\pm\pi$.
Since the signed areas are equal, the angles must be equal, and it
follows that the quadrilateral is cyclic, and thus a Darboux butterfly.

Note that if the equality holds for one (non-zero) value of $\ell$ then it holds for all values of $\ell$.

Finally, consider a non-planar quadrilateral $ABDC$ with the trivial monodromy (for some value of $\ell$). Assume that the monodromy along $ABD$ and $ACD$ are equal. Denote this monodromy by $M$. Then $M$ preserves the segments that lie in the plane $ABD$ and in the plane $ACD$, and hence, in their intersection, the line $AD$. In the plane $ABD$, the monodromy $M$ is given by  formula (\ref{matr}). If the horizontal axis is an eigendirection then $ab\sin(\alpha-\beta)=0$. This implies that the segments $AB$ and $BD$ are collinear, a contradiction.
\proofend

\begin{figure}[hbtp]
\centering
\includegraphics[height=1.2in]{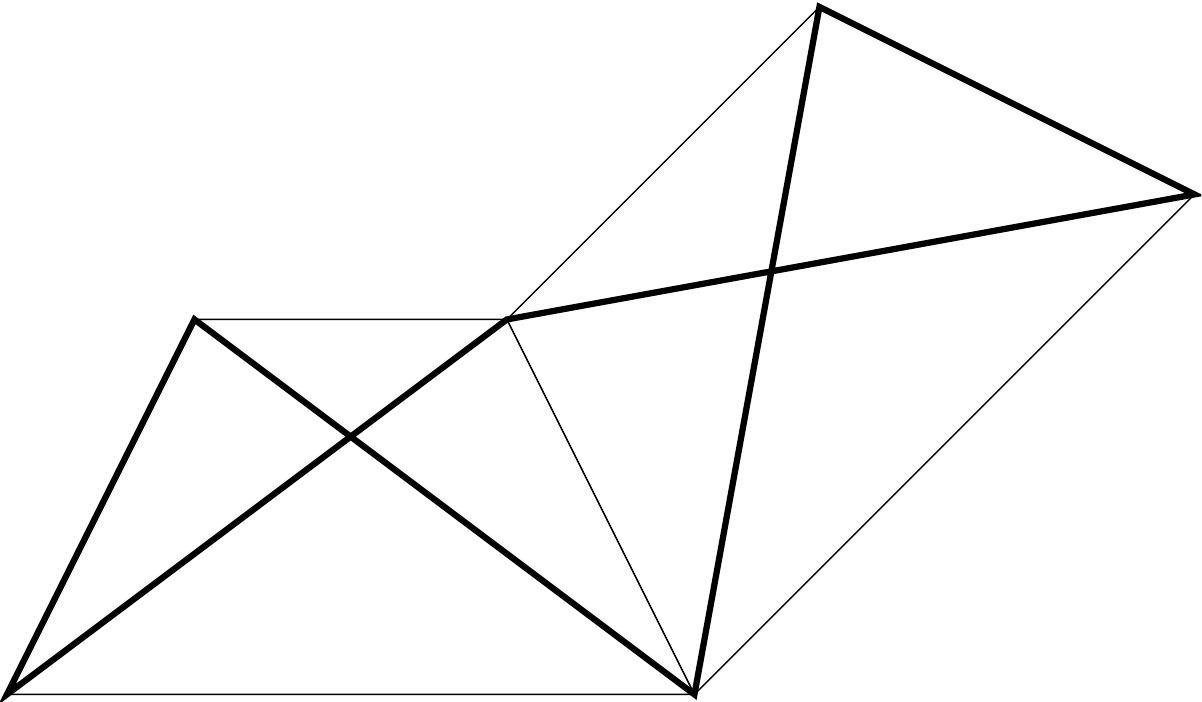}
\caption{Constructing polygons with identity monodromy}
\label{att}
\end{figure}

As a consequence of Butterfly Lemma, for every $n$, we can construct a family of $2n$-gons with identity monodromy for all values of $\ell$. These polygons are obtained by attaching Darboux butterflies to each other along the common sides, see Figure \ref{att}.

Now we are in a position to prove the rest of the theorems.
\medskip

\noindent{\bf Proof of Theorem \ref{conj}.} It follows from the Butterfly Lemma that, in Figure \ref{constr}, one has:
$$
M_{W_1 W_2,\lambda} = M_{V_2 W_2,\lambda} M_{V_1 V_2,\lambda}  M^{-1}_{V_1 W_1,\lambda}. 
$$ 
Taking the composition over the closed polygon $V$ yields the result.
\proofend
\medskip

\noindent{\bf Proof of Theorem \ref{comm}.} Consider the points $V_1, W_1, S_1$, and let $T_1$ be the point such that $V_1 W_1 T_1 S_1$ is a Darboux butterfly. Consider the discrete bicycle transformation of the segment $V_1 V_2$ along the Darboux butterfly $V_1 W_1 T_1 S_1$. According to the Butterfly Lemma, the resulting quadrilateral, say, $Q$, is closed and, according to Theorem \ref{conj}, it has the trivial monodromy (for any length parameter). Hence, by the Butterfly Lemma again, $Q$ is a Darboux butterfly as well. 

It is clear from Figure \ref{constr} that the discrete bicycle transformation of the segment $V_1 W_1$ along $V_1 V_2$ is the same as the discrete bicycle transformation of the segment $V_1 V_2$ along $V_1 W_1$. It follows that three of the vertices of $Q$ are $V_2, W_2$ and $S_2$. Denote the fourth vertex by $T_2$. 

A continuation of this process yields a closed polygon $T$ satisfying the assertion of the theorem.
\proofend
\medskip

\noindent{\bf Proof of Theorem \ref{rec}.} 
An equivalent description of recutting $V_i \mapsto V_i'$ is that the quadrilateral $V_{i-1} V_i V_{i+1} V_i'$ is a Darboux butterfly. 

To prove the first statement, we use Butterfly Lemma:
$$
M_{V_{i-1}V_i'V_{i+1},\lambda}=M_{V_{i-1}V_iV_{i+1}V_i'V_{i+1},\lambda}=M_{V_{i-1}V_iV_{i+1},\lambda}.
$$

For the second statement, let $W$ be be a polygon in the discrete bicycle correspondence with $V$. Let $V_i'W_i'$ be the discrete bicycle transformation of the segment $V_{i-1} W_{i-1}$ along the segment $V_{i-1} V_i'$. Since $V_{i-1} V_i V_{i+1} V_i'$ is a Darboux butterfly, the discrete bicycle  transformation takes $V_i'W_i'$ to $V_{i+1} W_{i+1}$. Thus the polygon $\dots W_{i-1} W_i' W_{i+1}\dots$ is in the discrete bicycle correspondence with $\dots V_{i-1} V_i' V_{i+1}\dots$

We want to show that the recutting of $W$ on $i$th vertex yields $W_i'$ or, equivalently, that $W_{i-1} W_i W_{i+1} W_i'$ is a Darboux butterfly. According to Butterfly Lemma, we need to show that the monodromy along the closed polygon $W_{i-1} W_i W_{i+1} W_i' $ is the identity.
Indeed, using that the monodromy of each Darboux butterfly is trivial, we obtain:
\begin{equation*}
\begin{split}
M_{W_{i-1} W_i W_{i+1} W_i' W_{i-1},\lambda} = M_{W_{i-1} V_{i-1} V_i W_i V_i V_{i+1} W_{i+1} V_{i+1} V_i' W_i' V_i' V_{i-1} W_{i-1},\lambda} \\=
 M_{W_{i-1} V_{i-1} V_i  V_{i+1}  V_i'  V_{i-1} W_{i-1},\lambda} =Id,
\end{split}
\end{equation*}
and we are done.
\proofend

\section{Integrals} \label{int}

As we mentioned earlier, the discrete bicycle transformation preserves the conjugacy equivalence class of the monodromy $M_{\lambda}$, thus yielding the monodromy integrals. These integrals do not change if a polygon is acted upon by an isometry of the ambient space. We plan to study the monodromy integrals in a forthcoming paper. In this section, we study the integrals introduced in \cite{Ad1,Ad2} as integrals of the recutting. One of these integrals, $J(V)$, is not preserved by isometries. The other integral, $A(V)$, was described, in the 3-dimensional case, in \cite{PSW}. 

Given a closed polygon $V$, consider the vector $J$ and the bivector $A$ given by the formulas
\begin{equation} \label{defAJ}
\begin{split}
&J(V)=\sum_{i} (|V_{i+1}|^2- |V_{i-1}|^2) V_{i}=\sum_i |V_i|^2 (V_{i-1}-V_{i+1}),\\
 &A(V)=\sum_{i} V_i \wedge V_{i+1},
 \end{split}
\end{equation}
where the sums are cyclic. 
In dimension 2, $A(V)$ is the signed area of the polygon $V$.

\begin{theorem} \label{AJ}
Both $A$ and $J$ are integrals of the discrete bicycle transformation.
\end{theorem}

As a preparation to the proof, we describe a discrete counterpart to the rear bicycle track (the middle curve with cusps in Figure \ref{Darboux}). 

\begin{figure}[hbtp]
\centering
\includegraphics[width=4in]{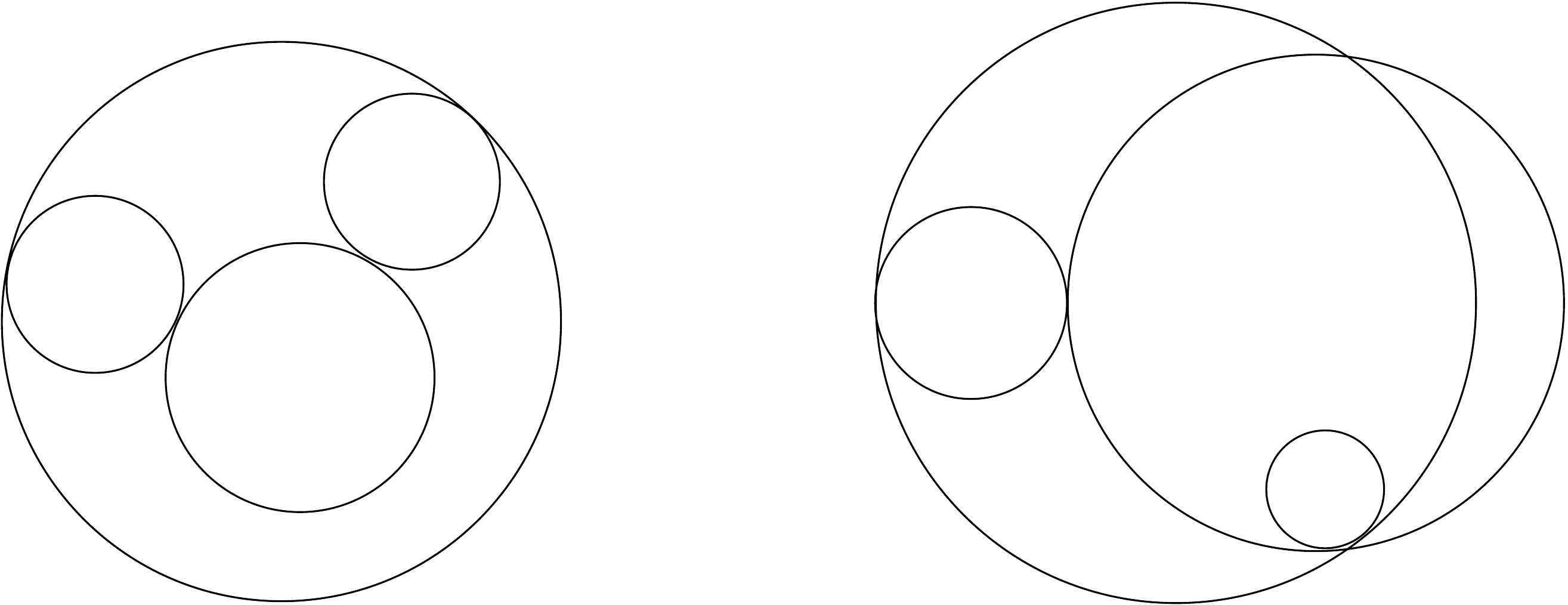}
\caption{An oriented and a non-oriented chain of four circles}
\label{circles}
\end{figure}

We shall consider collections of spheres such that the first one is tangent to the second, the second to the third,  ... , and the last one is tangent to the first. We call such a collection a {\it chain}. The radii of the spheres are signed. By convention, if two spheres have an exterior tangency then their radii have the same sign, and if the tangency is interior then the radii have the opposite signs. We allow infinite radii, that is, we consider hyperplanes as spheres as well. An infinite radius has no sign (equivalently, one may consider the curvatures, not excluding zero curvature form consideration). A chain is called {\it oriented} if one can choose the signs of the radii consistent with the sign convention. That is, a chain is oriented if and only if the number of interior tangencies is even, see Figure \ref{circles}. 

In what follows, we use half-integers as the indices for the centers of the spheres and of their radii. Consider an oriented chain of spheres with centers $P_j$ and signed radii $r_j$. Denote by $Q_i$ the tangency point of the adjacent spheres with centers $P_{i-\frac{1}{2}}$ and $P_{i+\frac{1}{2}}$. Let $V_i$ and $W_i$ be the two points on the line $P_{i-\frac{1}{2}} P_{i+\frac{1}{2}}$ located at distance $\ell$ from $Q_i$. The choice of labels is consistent for all $i$: if the segments $V_i W_i$ and $Q_i P_{i+\frac{1}{2}}$ have the same orientations then the segments $V_{i+1} W_{i+1}$ and $Q_{i+1} P_{i+\frac{1}{2}}$ have the opposite orientations, and vice versa.

\begin{lemma} \label{chainlemma}
The polygons $V$ and $W$ are in the discrete bicycle correspondence. Conversely, given polygons $V$ and $W$ in the discrete bicycle correspondence, let $P_{i+\frac{1}{2}}$ be the intersection point of the lines $V_{i+1} W_{i+1}$ and $V_i W_i$. Then there exists an oriented chain of spheres centered at points $P_j$, such that the tangency points $Q_i$ are the midpoints of the segments $V_i W_i$.
\end{lemma}

\begin{figure}[hbtp]
\centering
\includegraphics[width=3.5in]{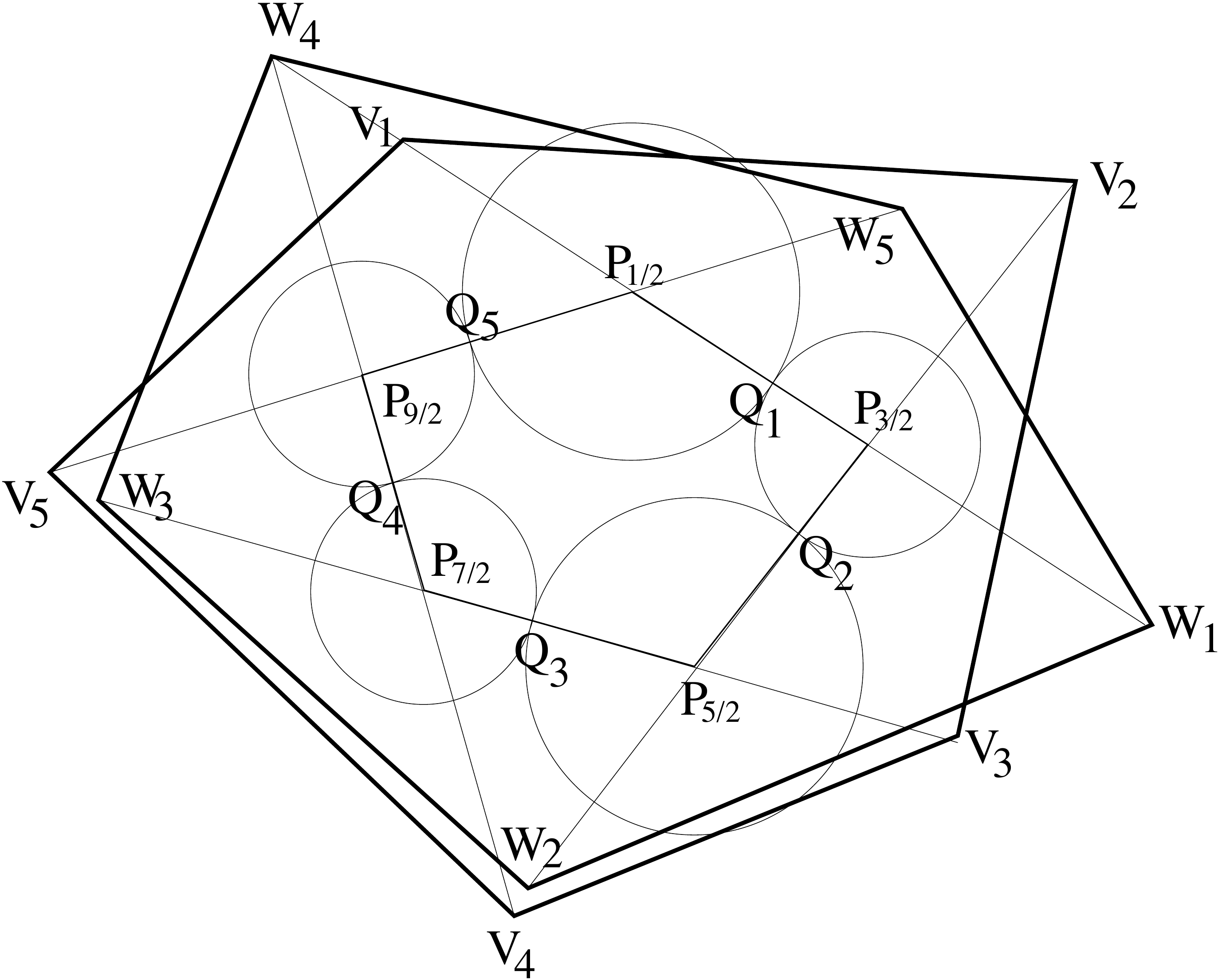}
\caption{Polygons $V$ and $W$ are in the discrete bicycle correspondence}
\label{chain}
\end{figure}

The construction is illustrated in Figure \ref{chain}.

\proof By construction, a homothety centered at $P_{i+\frac{1}{2}}$ takes $V_i$ to $W_i$ and $W_{i+1}$ to $V_{i+1}$. For example, in Figure \ref{chain}, the homothety with the coefficient 
$$
-\frac{\ell - r_{\frac{3}{2}}}{\ell + r_{\frac{3}{2}}},
$$
 centered at $P_{\frac{3}{2}}$, takes $V_1 W_2$ to $W_1 V_2$. Since $|V_i W_i| = |V_{i+1} W_{i+1}|=2\ell$, the quadrilateral $V_i W_i V_{i+1} W_{i+1}$ is a Darboux butterfly. 

Conversely, by construction, 
$$
|P_{i+\frac{1}{2}} W_i|= |P_{i+\frac{1}{2}} V_{i+1}|,\ \  |P_{i+\frac{1}{2}} V_i|= |P_{i+\frac{1}{2}} W_{i+1}|,
$$
hence $|P_{i+\frac{1}{2}} Q_i| = |P_{i+\frac{1}{2}} Q_{i+1}|:= r_{i+\frac{1}{2}}$ where $Q_i$ is the midpoint of the segment $V_i W_i$. The sphere with this radius passes through points $Q_i$ and $Q_{i+1}$ and is orthogonal to the lines $P_{i+\frac{1}{2}} Q_i$ and $P_{i+\frac{1}{2}} Q_{i+1}$. Thus one obtains a chain of spheres, and this chain is oriented.
\proofend

The polygon $Q$ is the discrete rear bicycle track. We apply  Lemma \ref{chainlemma} to prove Theorem \ref{AJ}.
\medskip

\noindent{\bf Proof of Theorem \ref{AJ}.} Given an oriented chain with centers at points $P_j$ and signed radii $r_j$ (where $j$ is half-integer), the tangency points $Q_i$ have the following coordinates:
$$
Q_{i}=\frac{r_{i+\frac{1}{2}} P_{i-\frac{1}{2}}+r_{i-\frac{1}{2}} P_{i+\frac{1}{2}}}{r_{i-\frac{1}{2}} +r_{i+\frac{1}{2}}}
$$
(note that this formula does not change if all radii are negated). Then the points $V_i$ and $W_i$ are given by the formula 
\begin{equation} \label{VW}
\frac{(r_{i+\frac{1}{2}}-\ell) P_{i-\frac{1}{2}}+(r_{i-\frac{1}{2}}+\ell) P_{i+\frac{1}{2}}}{r_{i-\frac{1}{2}} +r_{i+\frac{1}{2}}},
\end{equation}
where the positive $\ell$ gives $V_i$ and the negative $\ell$ gives $W_i$. 

To prove the invariance of $A$, we need to show that $A$ is an even function of $\ell$. Indeed, using formula (\ref{VW}), we find that the odd (linear in $\ell$) part of $A$ is
\begin{equation*}
\begin{split}
\sum \frac{(P_{i+\frac{1}{2}}-P_{i-\frac{1}{2}})\times (r_{i+\frac{3}{2}}P_{i+\frac{1}{2}}+r_{i+\frac{1}{2}} P_{i+\frac{3}{2}})}{(r_{i-\frac{1}{2}}+r_{i+\frac{1}{2}})(r_{i+\frac{1}{2}}+r_{i+\frac{3}{2}})}+\\ 
\frac{(r_{i+\frac{1}{2}}P_{i-\frac{1}{2}}+r_{i-\frac{1}{2}} P_{i+\frac{1}{2}})\times (P_{i+\frac{3}{2}}-P_{i+\frac{1}{2}})}{(r_{i-\frac{1}{2}}+r_{i+\frac{1}{2}})(r_{i+\frac{1}{2}}+r_{i+\frac{3}{2}})}=\\
\sum \frac{P_{i+\frac{1}{2}}\times P_{i+\frac{3}{2}}}{r_{i+\frac{1}{2}}+r_{i+\frac{3}{2}}} - \sum \frac{P_{i-\frac{1}{2}}\times P_{i+\frac{1}{2}}}{r_{i-\frac{1}{2}}+r_{i+\frac{1}{2}}} =0,
\end{split}
\end{equation*}
as needed.

To prove that $J$ is invariant, one  makes a similar computation. Let $e_i$ be the unit vector from $Q_i$ to $P_{i+\frac{1}{2}}$. 

One has $V_{i}=Q_{i}+\ell e_{i}$ and $W_{i}=Q_{i}-\ell e_{i}$. Hence 
$$
|V_{i}|^2 = \ell^2 + 2\ell\ Q_{i} \cdot e_{i} + |Q_{i}|^2.
$$
It follows that 
$$
|V_{i+1}|^2- |V_{i-1}|^2 = |Q_{i+1}|^2- |Q_{i-1}|^2 +2\ell\ (Q_{i+1}\cdot e_{i+1} - Q_{i-1}\cdot e_{i-1}),
$$
and the odd (linear in $\ell$) part of $J$ is
\begin{equation} \label{J1} 
\sum (|Q_{i+1}|^2- |Q_{i-1}|^2)\ e_{i} + 2 (Q_{i+1}\cdot e_{i+1} - Q_{i-1}\cdot e_{i-1})\ Q_{i}.
\end{equation}
Rewrite negative (\ref{J1}) as
\begin{equation} \label{J2}
\begin{split}
\sum |Q_{i}|^2 (e_{i+1}-e_{i-1}) + 2 Q_{i}\cdot e_{i} (Q_{i+1}-Q_{i-1})=\\
\sum |Q_{i}|^2 ((e_{i+1}+e_{i})-(e_{i}+e_{i-1})) + 2 Q_{i}\cdot e_{i} (Q_{i+1}-Q_{i-1}).
\end{split}
\end{equation}
Using the formulas
\begin{equation*}
\begin{split}
Q_{i-1}=Q_{i}-r_{i-\frac{1}{2}} (e_{i}+e_{i-1}),\ Q_{i+1}=Q_{i}+r_{i+\frac{1}{2}} (e_{i}+e_{i+1}),\\
 P_{i-\frac{1}{2}}=Q_{i}-r_{i-\frac{1}{2}}e_{i},\ P_{i+\frac{1}{2}}=Q_{i}+r_{i+\frac{1}{2}}e_{i},
\end{split}
\end{equation*}
rewrite (\ref{J2}) as
\begin{equation*}
\begin{split}
\sum (|Q_{i}|^2+2r_{i+\frac{1}{2}} Q_{i}\cdot e_{i}) (e_{i+1}+e_{i}) - (|Q_{i}|^2-2r_{i-\frac{1}{2}} Q_{i}\cdot e_{i}) (e_{i-1}+e_{i})=\\
\sum (|P_{i+\frac{1}{2}}|^2-r_{i+\frac{1}{2}}^2) (e_{i+1}+e_{i}) - \sum (|P_{i-\frac{1}{2}}|^2-r_{i-\frac{1}{2}}^2) (e_{i-1}+e_{i})  = 0,
\end{split}
\end{equation*}
as needed. 
\proofend

\begin{remark} \label{rels}
{\rm
One has the following relation between the integrals $A$ and $J$:
\begin{equation} \label{dirder}
D_{\xi} (J) (V) = -2 A(V) \cdot \xi = 2 \sum_i (V_i \cdot \xi)\ (V_{i-1}-V_{i+1}),
\end{equation}
where $D_{\xi}$ is the directional derivative along a vector $\xi$ and where dot is the Euclidean pairing of 2-vectors and vectors. Of course, (\ref{dirder}) is also an integral for every vector $\xi$.
}
\end{remark}

\begin{remark} \label{ccm}
{\rm
The integral $A$ is invariant under parallel translations, but $J$ is neither  invariant under parallel translations nor commutes with them. In dimension two, we adjust the integral $J$ so that it commutes with parallel translations and thus becomes a ``center", associated with a polygon. Namely, rotate $J(V)$ through $90^{\circ}$ and divide by four times the area: 
$$
\frac{1}{4A(V)} \big(\sum (y_i^2y_{i+1}-y_iy_{i+1}^2+x_i^2y_{i+1}-x_{i+1}^2y_i),\\ 
\sum (x_ix_{i+1}^2-x_i^2x_{i+1}+x_iy_{i+1}^2-x_{i+1}y_i^2)\big),
$$
where $V_i=(x_i,y_i)$ and the sums are cyclic. We call this point the {\it circumcenter of mass} of the polygon $V$ and denote it by $CCM(V)$. 

A justification of this terminology is as
follows. Consider a triangulation of the polygon $V$, and let $O_i$ be the circumcenter of $i$th triangle. Then $CCM(V)$ is the center of mass of the points $O_i$, taken with the weight equal to the (oriented) area of $i$th triangle. The result does not depend on triangulation. This construction is  mentioned in \cite{Ad1}; we plan to study it in detail in a forthcoming paper \cite{TT}.  

Let us mention, without proof, two properties of $CCM(V)$. First, if $V$ is an equilateral polygon then the circumcenter of mass coincides with the center of mass. This agrees with the observation, made in \cite{BMO2} that, in our terminology, the discrete bicycle transformation of an equilateral polygon preserves its center of mass. 

Second, in the continuous limit, as $V$ becomes a  curve $\gamma$, the circumcenter of mass of $V$ tends to the center of mass of the homogeneous lamina bounded by $\gamma$. As a consequence, the continuous bicycle transformation preserves the center of mass.
}
\end{remark}

We plan to study the monodromy integrals in a separate paper. We  comment on these integrals in dimension two in the next section.


\section{In the plane} \label{plane}

\begin{figure}[hbtp]
\centering
\includegraphics[width=2.3in]{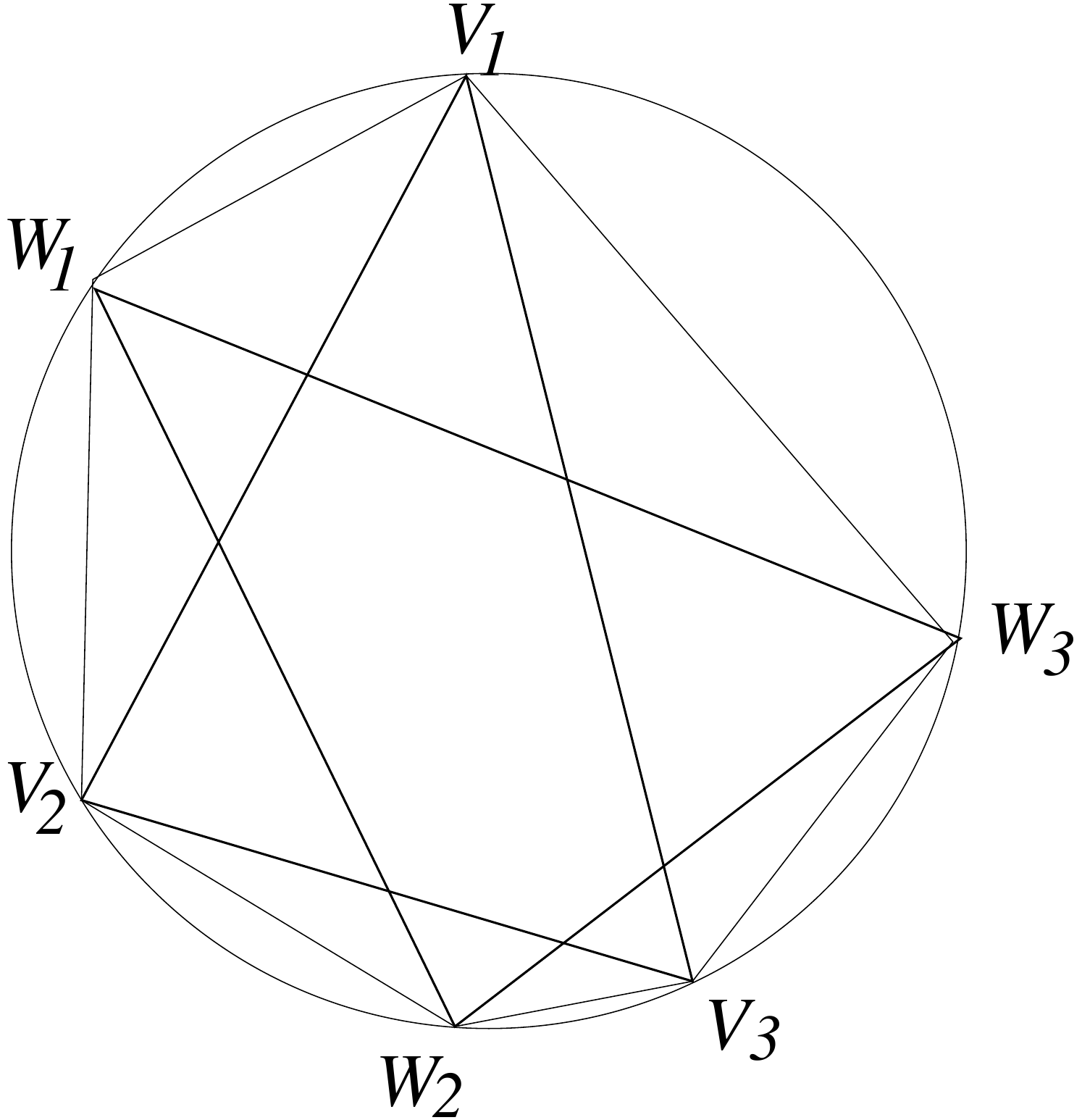}
\caption{Triangles $V$ and $W$ are in the discrete bicycle correspondence}
\label{rot}
\end{figure}

In this section, we consider the discrete bicycle transformation in the plane.
We start with a simple observation: for an inscribed polygon, a rotation about the circumcenter is a discrete bicycle transformation, see Figure \ref{rot}. 

Our first result concerns convex inscribed polygons.

\begin{theorem} \label{inscr}
Let $V$ be a convex inscribed polygon, and let $d$ be the diameter of the circumcircle. The discrete bicycle monodromy $M_{V,\ell}$ is elliptic for $\ell>d$, parabolic for $\ell=d$, and hyperbolic for $\ell \in (0,d)$. In the last case, the discrete bicycle transformation is a rotation about the circumcenter.
\end{theorem}

\proof As we mentioned, if $\ell \in (0,d)$ then a rotation about the circumcenter is a discrete bicycle transformation. 

Let $W=T_{\ell}(V)$. Then $W$ has the same perimeter and the same oriented area as $V$, see Theorem \ref{AJ}. It is known that, among polygons with given side lengths, there exists a unique area maximizing one, and this is an inscribed convex polygon. It follows that $W$ is inscribed, and hence congruent to $W$. It follows from Theorem \ref{AJ} that the circumcenter of $W$ coincides with that of $V$, see Remark \ref{ccm}. It follows that $W$ is a rotation of $V$ about the circumcenter.
\proofend

In particular, Theorem \ref{inscr} completely described the discrete bicycle transformation on triangles. 

Next we consider a $2k$-gon whose sides lie, in an alternating fashion, on two concentric circles. In the limiting  case, the two concentric circles may become two parallel lines. 

\begin{proposition} \label{concentric}
Let $C_1$ and $C_2$ be concentric circles with the center $O$ (or parallel lines). Let the odd vertices of a $2k$-gon lie on $C_1$ and the even ones on $C_2$. Let $W_1$ be a point of $C_2$. Then the discrete bicycle transformation of $V$ with the initial segment $V_1 W_1$ is a closed $2k$-gon whose odd vertices lie on $C_2$ and the even ones on $C_1$. The second iteration of this discrete bicycle transformation sends $V$ to an isometric polygon.
\end{proposition} 

\begin{figure}[hbtp]
\centering
\includegraphics[width=1.8in]{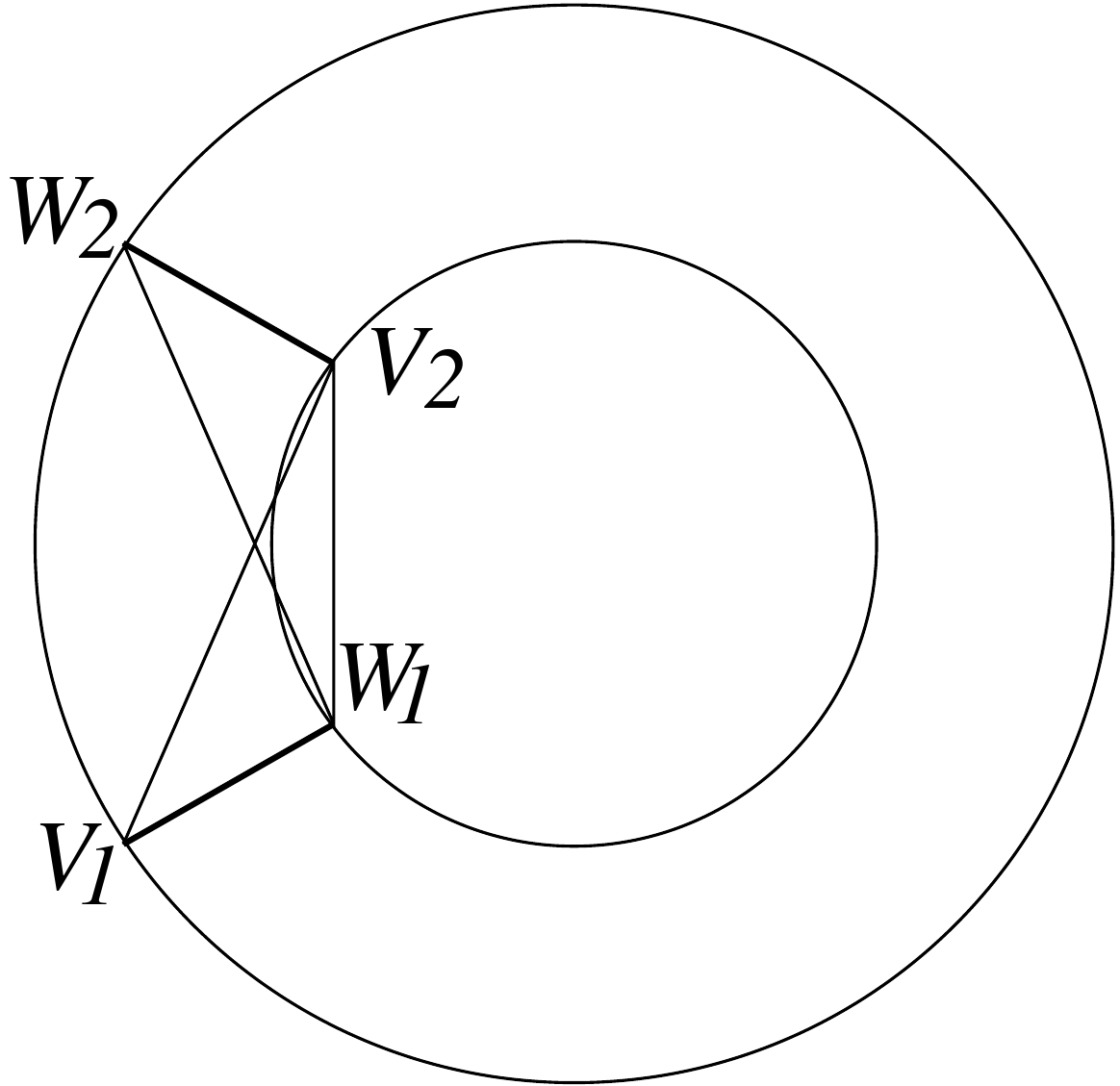}
\caption{$W_2$ is the reflection of $V_1$ in the perpendicular bisector of  $W_1 V_2$}
\label{conc}
\end{figure}

\proof Reflect $V_1$ in the perpendicular bisector of the segment $W_1 V_2$ to obtain $W_2$, and continue in the same way, see Figure \ref{conc}. Let the lower case letters denote the angular coordinates of the respective points. Then 
$$
w_2=w_1+v_2-v_1,\ w_3=w_2+ v_3-v_2=w_1 + v_3-v_1, 
$$
etc. It follows that $w_{2k+1} = w_1+ v_{2k+1}-v_1=w_1$, hence the polygon $W$ is closed.

We see that the discrete bicycle transformation ${\cal T}$ is the composition of two commuting transformations: the rotation through the angle $w_1-v_1$, and the involution that interchanges the points of $C_1$ and $C_2$ on the same radial ray. Hence ${\cal T}^2$ is a rotation.

The argument for parallel lines is analogous, and the resulting polygon $W$ is obtained from $V$ by a glide reflection. In this case, the orbit of the polygon is unbounded.
\proofend

Note that, in this construction, the polygon $Q$, whose vertices are the midpoints of the segments $V_i W_i$ (see  Lemma \ref{chainlemma}), is inscribed in a circle with the center $O$. We also have the following consequence of the proof.

\begin{corollary} \label{rhombus}
If the polygon in the preceding Proposition is a rhombus then its image under the bicycle transformation is a congruent rhombus.
\end{corollary}

Now we discuss the monodromy integrals for plane polygons. The monodromy along a side is given by formula (\ref{mono}); the full monodromy $M$ is the product of these monodromies over the consecutive sides of the polygon. The monodromy is defined only up to a multiplicative factor, and the invariant quantity is 
$$
\frac{\mbox{Tr}^{2}(M)}{\det(M)},
$$ 
considered as a function of $\ell$. 
Note that the determinant of the matrix (\ref{mono}) equals $\ell^2-a^2$, that is, is also an integral. Thus $\mbox{Tr}(M)$ is an integral.

\begin{proposition} \label{intpol}
Consider a  $k$-gon whose sides have the lengths $a_1,\dots,a_k$ and the directions $\alpha_1,\dots,\alpha_k$.
Then
$$
{\rm Tr} (M)=2(\ell^k + c_1 \ell^{k-1}+c_2\ell^{k-2}+\dots + c_k)
$$
with all odd coefficients $c_1, c_3, \dots$ equal to zero. If $k$ is even then the free term $c_k$ equals 
$$
a_1 \ldots  a_k \cos(\alpha_{1}-\alpha_{2}+\dots -\alpha_{k}).
$$
One also has:
$$
c_2=- \frac{1}{2}\sum a_i^2.
$$
\end{proposition}

\proof One has
$$
M=\prod_{i=1}^k (\ell E + a_i A(\alpha_i))
$$
where 
$$
A(\alpha)=
\left(\begin{array}{cc}
\cos\alpha & -\sin\alpha\\
-\sin\alpha & -\cos\alpha
\end{array}\right).
$$
Therefore 
$$
{\rm Tr} (M) = \sum_{j=0}^k \ell^{k-j} a_{i_1}\dots a_{i_j} {\rm Tr} (A(\alpha_{i_1})\dots A(\alpha_{i_j})).
$$

Notice that 
\begin{equation} \label{prod}
A(\alpha) A(\beta) =
\left(\begin{array}{cc}
\cos(\alpha-\beta) & \sin(\alpha-\beta)\\
-\sin(\alpha-\beta) & \cos(\alpha-\beta)
\end{array}\right),
\end{equation}
a rotation matrix. More generally, the product of an odd number of the matrices $A(\alpha_i)$ is traceless, and the product of an even number is a rotation through the alternating sum of the respective angles. This implies the first two claims.

For the last claim, let $u_1,\dots, u_k$ be the vectors of the sides of the polygon. Using (\ref{prod}), we find that
$$
c_2=\sum_{i<j} u_i \cdot u_j.
$$
One has: $\sum u_i =0$. Taking dot with itself yields:
$$
0=\sum u_i \cdot u_i + 2 \sum_{i<j} u_i \cdot u_j.
$$
Thus 
$$
c_2= -\frac{1}{2}\sum a_i^2,
$$
as claimed.
\proofend

\begin{corollary} \label{alt}
The quantity
$\cos(\alpha_{1}-\alpha_{2}+\dots -\alpha_{k})$
is an integral of the discrete bicycle transformation on even-gons.
\end{corollary}

\begin{figure}[hbtp]
\centering
\includegraphics[width=1.7in]{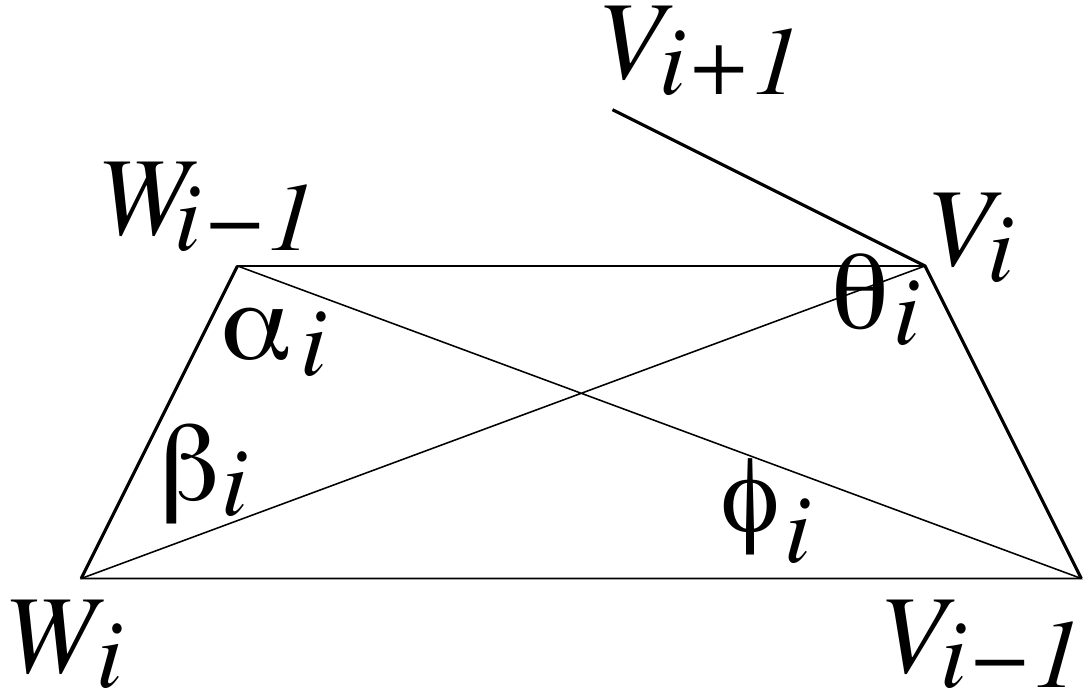}
\caption{Notations for Proposition \ref{angrec}}
\label{angles}
\end{figure}

Let polygons $V$ and $W$ be in the discrete bicycle correspondence. Let 
$$
\alpha_i = \angle V_{i-1} V_i W_i = \angle V_{i-1} W_{i-1} W_i,
$$
see Figure \ref{angles}. If one knows the cyclic sequence of angles $\alpha_i$ then one can construct $W$ from $V$: indeed, the lengths of all the segments $V_iW_i$ are equal to $2\ell$.

The angles $\alpha_i$ satisfy a first order nonlinear difference equation with periodic coefficients. Let $\theta_i = \angle V_{i-1}V_iV_{i+1}$ and  $c_i=|V_{i-1}V_i|$.

\begin{proposition} \label{angrec}
One has
\begin{equation} \label{difference}
2\ell \cos \left(\frac{\alpha_i-\alpha_{i-1}+\theta_{i-1}}{2}\right) =c_i \cos \left(\frac{\alpha_i+\alpha_{i-1}-\theta_{i-1}}{2}\right).
\end{equation}
\end{proposition}

\proof Let 
$$
\beta_i=\angle W_{i-1} V_{i-1} V_i = \angle W_{i-1} W_i V_i,\ \ \phi_i = \angle W_{i-1} V_{i-1} W_i = \angle V_i W_i V_{i-1}.
$$
Then $2\phi_i=\pi-\alpha_i-\beta_i$. Since $\angle W_i V_i V_{i+1} = \beta_{i+1}$, one  has $\beta_{i+1}=\theta_{i} - \alpha_i$. Therefore
\begin{equation} \label{angeq}
\phi_i=\frac{\pi}{2} - \frac{\alpha_i-\alpha_{i-1}+\theta_{i-1}}{2},\ \beta_i+\phi_i=\frac{\pi}{2} - \frac{\alpha_{i-1}+\alpha_i-\theta_{i-1}}{2}.
\end{equation}

By Sine Rule in triangle $V_{i-1} V_i W_i$,
$$
\frac{2\ell}{\sin(\beta_i+\phi_i)} = \frac{c_i}{\sin \phi_i}, 
$$
or
$$
2\ell \cos \left(\frac{\alpha_i-\alpha_{i-1}+\theta_{i-1}}{2}\right) =c_i \cos \left(\frac{\alpha_i+\alpha_{i-1}-\theta_{i-1}}{2}\right),
$$
as claimed. 
\proofend

As an application of Proposition \ref{angrec}, we compute the eigenvalue of the fixed point of the monodromy map of the polygon $V$ corresponding to the pair of polygons $V,W$ in the discrete bicycle correspondence. Since  the monodromy is a M\"obius transformation, the eigenvalues of its two fixed points are reciprocals of each other. 

\begin{theorem} \label{eigen}
The eigenvalue in question equals
$$
\prod_{i=1}^n \frac{|V_{i-1} W_{i}|}{|V_{i} W_{i-1}|} = \prod_{j=1/2}^{n+1/2} \frac{|\ell+r_j|}{|\ell-r_j|}.
$$
In particular, the monodromy is parabolic if and only if 
$$
\prod_{i=1}^n |V_{i-1} W_{i}| = \prod_{i=1}^n |V_{i} W_{i-1}|
\quad {\rm or} \quad
\prod_{j=1/2}^{n+1/2} |\ell+r_j| = \prod_{j=1/2}^{n+1/2} |\ell-r_j|.
$$

\end{theorem}

\proof To compute the eigenvalue, one linearizes equation (\ref{difference}): if $u_i$ is a variation of $\alpha_i$ then the linearization is as follows:
$$
2\ell (u_i-u_{i-1}) \sin\left(\frac{\alpha_i-\alpha_{i-1}+\theta_{i-1}}{2}\right) = c_i (u_i+u_{i-1}) \sin\left(\frac{\alpha_i+\alpha_{i-1}-\theta_{i-1}}{2}\right),
$$
and hence
\begin{equation*}
\begin{split}
u_i \left[2\ell \sin\left(\frac{\alpha_i-\alpha_{i-1}+\theta_{i-1}}{2}\right) - c_i \sin\left(\frac{\alpha_i+\alpha_{i-1}-\theta_{i-1}}{2}\right)\right] = \\
u_{i-1}  \left[2\ell \sin\left(\frac{\alpha_i-\alpha_{i-1}+\theta_{i-1}}{2}\right) + c_i \sin\left(\frac{\alpha_i+\alpha_{i-1}-\theta_{i-1}}{2}\right)\right].
\end{split}
\end{equation*}
By elementary geometry of the trapezoid in Figure \ref{angrec} and formulas (\ref{angeq}), one has:
\begin{equation*}
\begin{split}
2\ell \sin\left(\frac{\alpha_i-\alpha_{i-1}+\theta_{i-1}}{2}\right) = \frac{1}{2} (|V_{i-1} W_{i}|+|V_{i} W_{i-1}|),\\
c_i \sin\left(\frac{\alpha_i+\alpha_{i-1}-\theta_{i-1}}{2}\right) = \frac{1}{2} (|V_{i-1} W_{i}|-|V_{i} W_{i-1}|).
\end{split}
\end{equation*}
Therefore
$$
u_i |V_{i} W_{i-1}| = u_{i-1} |V_{i-1} W_{i}|,
$$
which implies the first formula for the eigenvalue. 

For the second formula, note that a homothety centered at point $P_{i+1/2}$ takes segment $V_i W_{i+1}$ to segment $V_{i+1} W_i$, see Figure \ref{chain}. The coefficient of this homothety is $|\ell+r_{i+1/2}|/|\ell-r_{i+1/2}|$, and we obtain the second formula for the eigenvalue.

It remains to notice that the monodromy is parabolic if and only if the two reciprocal  eigenvalues coincide.
\proofend

\begin{remark} \label{contcase}
{\rm The continuous analogs of Proposition \ref{angrec} and Theorem \ref{eigen} are contained in \cite{LT}. Namely, the continuos version of (\ref{difference}) is the differential equation
$$
\frac{d\alpha}{dx}+ \frac{\sin\alpha}{\ell}=\kappa(x)
$$
where $\alpha(x)$ is the angle made by the bicycle frame with the front wheel trajectory, $x$ is the arc length parameter along this trajectory, and $\kappa(x)$ is the curvature of this curve. The endpoint of the segment of length $\ell$ describes the rear wheel trajectory.

The continuos version of Theorem \ref{eigen} states that the eigenvalues of the bicycle monodromy are 
$
e^{\pm length(\gamma)}
$
where $\gamma$ is the rear wheel trajectory, and the length is algebraic: the sign changes after one traverses a cusp. In particular, the monodromy is parabolic if and only if the rear track has zero length. 
}
\end{remark}

\section{Case study: plane quadrilaterals} \label{quad}

In this section, we describe the dynamics of the discrete bicycle transformation on plane quadrilaterals.

\begin{figure}[hbtp]
\centering
\includegraphics[width=2.3in]{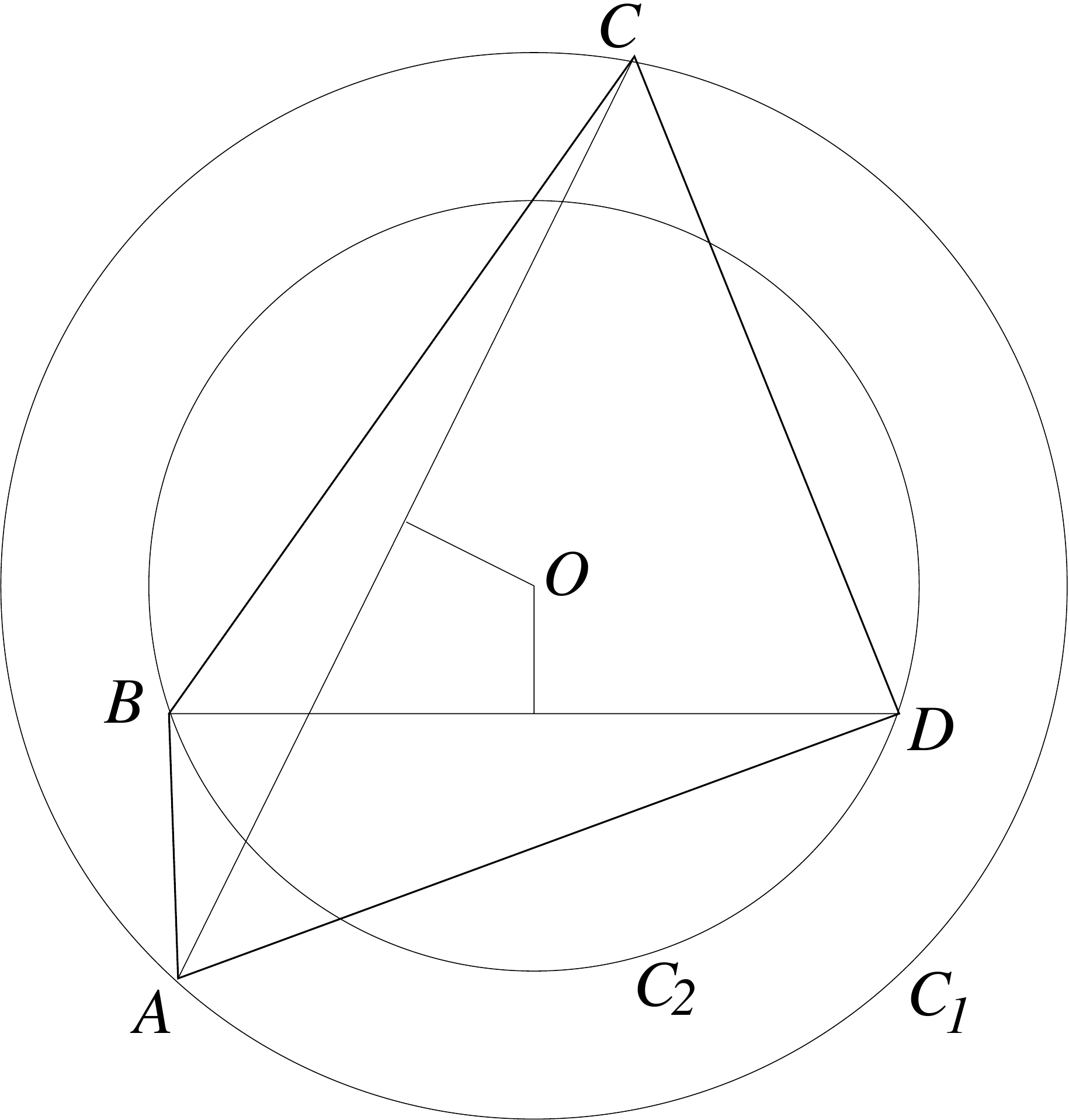}
\quad\quad
\includegraphics[width=2in]{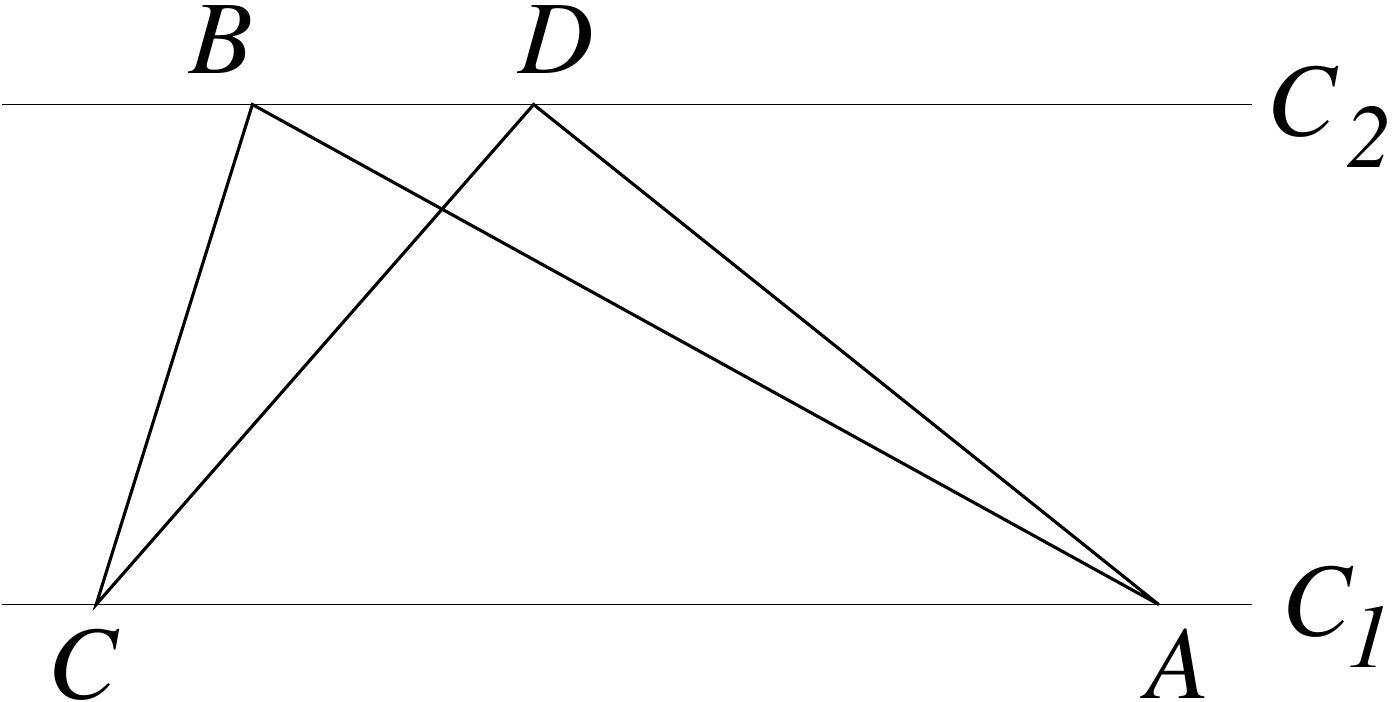}
\caption{Two types of quadrilaterals: the circumcenter is finite or infinite}
\label{quads}
\end{figure}

We have a trichotomy according to the position of the circumcenter of mass, see Remark \ref{ccm}. Consider a quadrilateral $ABCD$. The first case is when the diagonals $AC$ and $BD$ are not parallel. Let $O$ be the intersection point of the perpendicular bisectors of these diagonals, see Figure \ref{quads} on the left. 

\begin{lemma} \label{CCMass}
$O$ is the circumcenter of mass of the quadrilateral $ABCD$.
\end{lemma}

\proof The circumcenters of the triangles $ABD$ and $BCD$ lie on the perpendicular bisector of the segment $BD$, and the circumcenters of the triangles $ABC$ and $ACD$ lie on the perpendicular bisector of the segment $AC$. Hence $O=CCM(ABCD)$.
\proofend

In the first case, $A$ and $C$ lie on one circle, say, $C_1$, and $B$ and $D$ on another circle, $C_2$, centered at $O$. Denote their radii by $r_1$ and $r_2$, and assume that $r_1 \geq r_2$.

The second case is when the diagonals are parallel but the quadrilateral is not a Darboux butterfly, see Figure \ref{quads} on the right. In this case, the two concentric circles are replaced by two parallel lines, and the center $O$ is at infinity. Although both radii are infinite,  their difference $r_1-r_2$ is still  defined and equals the distance between the parallel lines. Note that, in this case, the quadrilateral $ABCD$ has zero area.

The third case is when the quadrilateral is a Darboux butterfly. In this case, there exists an infinite family of pairs of concentric circles $C_1, C_2$ such that $A,C \in C_1$ and $B,D \in C_2$. The centers of these circles lie on the common perpendicular bisector of the segments $AC$ and $BD$, including the point at infinity, when the circles become parallel lines.

\begin{theorem} \label{mainquad}
Let $ABCD$ be a quadrilateral. If $ABCD$ is not a Darboux butterfly then the discrete bicycle monodromy about the quadrilateral is elliptic for $\ell \in (0, r_1 - r_2) \cup (r_1+r_2, \infty)$, hyperbolic for $\ell \in (r_1-r_2,r_1+r_2)$, and parabolic for $\ell = r_1\pm r_2$. For $\ell$ in the hyperbolic or parabolic range, the discrete bicycle correspondence is induced by a point $A' \in C_2$, as described in Proposition \ref{concentric}. 
If $ABCD$ is a Darboux butterfly then the monodromy is the identity. For every starting point $A'$, there exists a circle (or straight line) $C_2$ that passes through $A'$, and the  discrete bicycle correspondence is again described by Proposition \ref{concentric}. 
\end{theorem}

\proof If $\ell \in [r_1-r_2,r_1+r_2]$ then there exist two points $A' \in C_2$ such that $|AA'|=\ell$ (these two points coincide for $\ell=r_1\pm r_2$), and Proposition \ref{concentric} describes the discrete bicycle transformation. 

Conversely, assume that $A'B'C'D'$ is a discrete bicycle transformation of $ABCD$. Let $l_1, l_2, l_3$ and $l_4$ be the perpendicular bisectors of the segments $A'B, B'C, C'D$ and $D'A$, respectively. Let $R_i$ be the reflection in the line $l_i,\ i=1,2,3,4$.
By definition of the bicycle monodromy, 
$$
B'=R_1 (A),\ C'= R_2(B),\ D' = R_3(C),\ A'=R_4(D),
$$
see Figure \ref{conc}. Note also that
$$
B=R_1 (A'),\ C= R_2(B'),\ D = R_3(C'),\ A=R_4(D').
$$

We claim that the lines $l_1,l_2,l_3,l_4$ are concurrent (as a particular case, the four lines may be parallel). 

Consider  the composition $F=R_3\circ R_2\circ R_1$: it is either a reflection or a glide reflection. We claim that the former is the case. Two given congruent line segments $AA'$,$D'D$ are related by just one odd isometry. Since $AA'D'D$ is an isosceles trapezoid, this isometry is a reflection.


Since $R_3\circ R_2\circ R_1$ is a reflection, the lines $l_1,l_2$ and $l_3$ are concurrent. Applying the same argument to $l_2,l_3, l_4$, we conclude that all four lines are concurrent.

To fix ideas, let us assume that the intersection point of the lines $l_1,l_2,l_3,l_4$ is finite (the case of parallel lines is similar). Denote this point by $Q$. We claim that $Q=O$, the circumcenter of the quadrilateral $ABCD$. 

Indeed, $R_2 \circ R_1 (A)=C$, hence $Q$ lies on the perpendicular bisector of the diagonal $AC$. Likewise, $R_3 \circ R_2 (B) =D$, hence $Q$ lies on the perpendicular bisector of the diagonal $BD$. Thus $Q=O$.

Since $A'=R_1(B)$, it follows that $A'\in C_2$, and we are in the situation of Proposition \ref{concentric}.

It remains to consider the case of a Darboux butterfly. For any starting point $A'$, we can find a circle $C_2$ through $A',B$ and $D$ with the center $O$ on the perpendicular bisector of the segments $AC$ and $BD$. Then another circle $C_1$, centered at $O$, passes through $A$ and $C$, and we are in the situation described in Proposition \ref{concentric}.
\proofend

\begin{remark} \label{alter}
{\rm The   preceding argument provides an alternative proof of the fact that the monodromy of a Darboux butterfly is the identity for all $\ell$.
}
\end{remark}

We now discuss an application of Theorem \ref{mainquad} to the following  problem in ``bicycle mathematics". Suppose one is given two closed curves, the front and rear bicycle tracks. Can one always determine in which direction the bicycle went? Usually, one can, but sometimes one cannot: consider, for example, two concentric circles. 

Describing such pairs of ``ambiguous" bicycle tracks is an interesting and difficult problem, and only partial results are available. This problem is equivalent to Ulam's problem of describing uniform (2-dimensional) bodies that float in equilibrium in all positions. We refer to \cite{BMO1,BMO2,Fi,We1,We2} for the literature on this intriguing topic.

A discrete version of this problem was introduced in \cite{Ta1}. Define a {\it bicycle $(n,k)$-gon} as an equilateral $n$-gons whose $k$-diagonals have equal length. More precisely, if the polygon is $V_1 V_2 \dots V_n$ then we require that $V_i V_{i+1} V_{i+k+1} V_{i+k}$ be a Darboux butterfly for all $i$ (as usual, the indices are understood cyclically). The problem is to describe bicycle $(n,k)$-gons, in particular, to determine for which pairs $(n,k)$ such a polygon must be regular. See also \cite{CC,Cs}.

\begin{figure}[hbtp]
\centering
\includegraphics[width=1.4in]{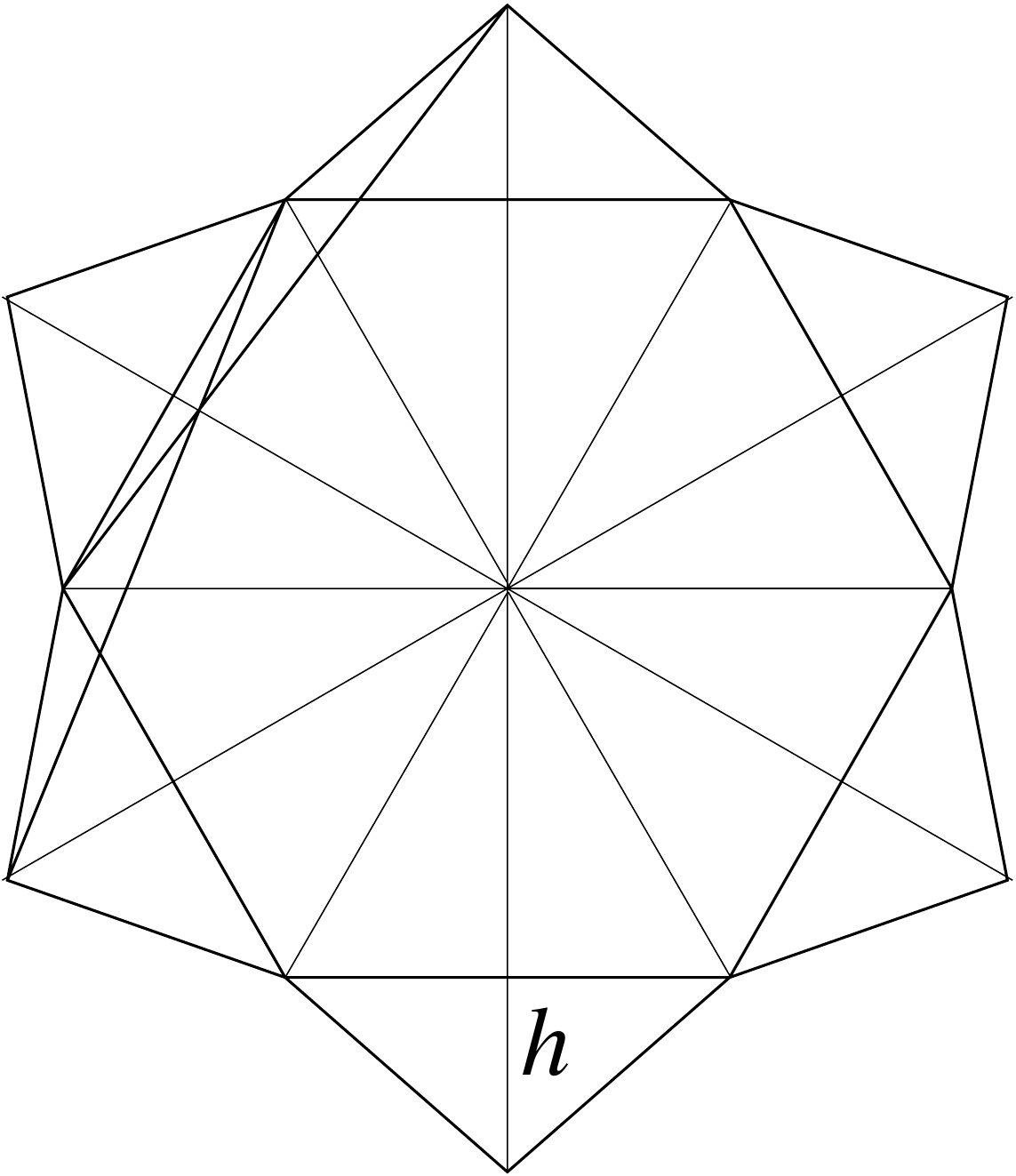}
\caption{A bicycle $(12,3)$-gon: $h$ is a parameter of the construction}
\label{bipol}
\end{figure}

For example, it is shown in \cite{Ta1} that bicycle $(n,2)$-gons, $(2k+1,k)$-gons, and $(3k,k)$-gons are regular. On the other hand, an example  of a non-regular bicycle polygon is shown in Figure \ref{bipol}. This construction generalizes to all pairs $(n,k)$ and yields 1-parameter families of bicycle $(n,k)$-gons with even $n$ and odd $k$. Note that the even and the odd vertices of a polygon in Figure \ref{bipol} lie on two concentric circles and that the polygons have dihedral symmetry.

Let $\ell$ be the length of the $k$-diagonal of a bicycle $(n,k)$-gon, and let $S$ be the cyclic relabeling of the vertices: $V_i \mapsto V_{i+1}$.
One can restate the definition in terms of the discrete bicycle transformation ${\cal T}_{\ell}$:  $V$ is a bicycle $(n,k)$-gon if ${\cal T}_{\ell} (V)=S^k(V)$.

The next result is a further step toward classification of bicycle polygons. 

\begin{theorem} \label{rigid}
If $k$ is even then a bicycle $(4k,k)$-gon is regular. If $k$ is odd then the even vertices of a bicycle $(4k,k)$-gon are equally spaced on a circle and its odd vertices are equally spaced on a concentric circle, that is, the polygon is obtained from a regular $2k$-gon by the construction depicted in Figure \ref{bipol}.
\end{theorem}

\proof
Given a bicycle $(4k,k)$-gon $V$, consider the rhombus $V_0 V_{k} V_{2k} V_{3k}$. The discrete bicycle transformation with the length parameter $V_0 V_1$ takes this rhombus to $V_1 V_{k+1} V_{2k+1} V_{3k+1}$, to $V_2 V_{k+2} V_{2k+2} V_{3k+2}$, and so on. 

Let $O$ be the center of the rhombus $V_0 V_{k} V_{2k} V_{3k}$, and let $C_1$ and $C_2$ be the concentric circles centered at $O$ such that $V_0, V_{2k} \in C_1$ and  $V_k, V_{3k} \in C_2$. By Corollary \ref{rhombus}, all the consecutive rhombi are congruent, and 
$ V_1 \in C_2, V_2 \in C_1, V_3 \in C_2, V_4 \in C_1$, etc. 

Therefore, if $k$ is even, then $V_k \in C_1$, and hence $C_1=C_2$. It follows that the rhombus is a square and $V$ is a regular $4k$-gon. If $k$ is odd then the even vertices of $V$ form a regular $2k$-gon inscribed into $C_1$, and the odd ones form a regular $2k$-gon inscribed into $C_2$. Thus $V$ is obtained from a regular $2k$-gon by the construction in Figure \ref{bipol}.  
\proofend

\bigskip
{\bf Acknowledgments}. We have discussed the discrete bicycle transformation with many a mathematician, and we are grateful to them all. In particular, it is a pleasure to acknowledge interesting discussions with I. Alevi, A. Bobenko,  T. Hoffmann,  U. Pinkall, B. Springborn, Yu. Suris, and  A. Veselov. This project originated during the program Summer@ICERM 2012; we are grateful to ICERM for support and hospitality. 
S. T. was partially supported by the NSF grant DMS-1105442.

\end{document}